\documentclass[11pt, final]{amsart}
\usepackage{amssymb}
\usepackage[mathscr]{eucal}

\theoremstyle{plain}
\newtheorem*{theorem}{Theorem}
\newtheorem*{corollary}{Corollary}

\theoremstyle{remark}
\newtheorem{remark}{Remark}

\def\cA{\mathscr{A}}
\def\C{\mathbb{C}}
\def\I{\mathscr{I}}
\def\l{\lambda}
\def\la{\langle}
\def\N{\mathbb N}
\def\ot{\otimes}
\def\ra{\rangle}

\begin{document}
\title[]{Multiplicative maps on ideals of operators which are local
automorphisms}
\author{LAJOS MOLN\'AR}
\address{Institute of Mathematics\\
         Lajos Kossuth University\\
         4010 Debrecen, P.O.Box 12, Hungary}
\email{molnarl@math.klte.hu}
\thanks{  This research was supported from the following sources:\\
          1) Hungarian National Foundation for Scientific Research
          (OTKA), Grant No. T--016846 F--019322,\\
          2) A grant from the Ministry of Education, Hungary, Reg.
          No. FKFP 0304/1997}
\begin{abstract}
We present the following reflexivity-like result concerning the
automorphism group of the $C^*$-algebra $B(H)$, $H$ being a separable
Hilbert
space. Let $\phi:B(H)\to B(H)$ be a multiplicative map (no linearity or
continuity is assumed) which can be approximated at every point by
automorphisms of $B(H)$ (these automorphisms may, of course, depend on
the point) in the operator norm. Then $\phi$ is an automorphism of the
algebra $B(H)$.
\end{abstract}
\maketitle

\section{Introduction}

Motivated by a problem of Larson
\cite[Some concluding remarks (5), p. 298]{Lar} and a result of
Kadison \cite{Kad}, in a series of papers we investigated the problem of
the reflexivity of the automorphism groups of various operator algebras
(see \cite{MolBatty}, \cite{MolStud}, \cite{MolJFA} and the references
therein as well). This concerns the following question. Let $\cA$ be  a
$C^*$-algebra (or, more generally, a Banach algebra).
Let $\phi :\cA \to \cA$ be any bounded linear transformation. Suppose
that for every $A\in \cA$ there exists a sequence $(\phi_n)$ of
automorphisms of $\cA$ such that $\phi(A)=\lim_n \phi_n(A)$. Does it
follow that $\phi$ is an automorphism?
If the answer is 'yes', then we say that the automorphism group of $\cA$
is topologically reflexive. If the condition that for every $A\in \cA$
there exists an automorphism $\phi_A(A)$ of $\cA$ such that
$\phi(A)=\phi_A(A)$ implies that $\phi$ is an automorphism, then the
automorphism group of $\cA$ is called algebraically reflexive.
In both cases, that is when speaking about topological or algebraic
reflexivity, we have assumed that the map $\phi$ is linear and tried
to conclude that then its multiplicativity and bijectivity follow. It
is now a natural question that what happens if we assume the
multiplicativity of $\phi$ and try to get its linearity and bijectivity.
The interesting problem that when does the mutiplicativity imply
additivity has been investigated in several papers. See, for example,
\cite{Hak}, \cite{Sem} which treat
the additivity of semigroup automorphisms (we emphasize that bijectivity
is assumed there) of operator algebras and \cite{Mar} for a general
algebraic result on surjective semigroup homomorphisms onto
prime rings. The aim of this note is to present the multiplicative
counterpart
of our recent topological reflexivity result for the automorphism group
of $B(H)$ obtained in \cite{MolBatty}.

Let us now fix the notation and the concepts that we shall use
throughout.
If $\cA$ is a *-algebra, then by an automorphism of $\cA$ we mean a
linear and multiplicative bijection of $\cA$ onto itself which preserves
the *-operation.
In what follows let $H$ be a separable Hilbert space and denote by
$B(H)$ the $C^*$-algebra of all bounded linear operators acting on $H$.
The ideal of finite rank operators and that of compact operators on $H$
are denoted by $F(H)$ and $C(H)$, respectively.
It is a well-known result due to Calkin that for every nontrivial ideal
$\I$ of $B(H)$ we have $F(H) \subset \I \subset C(H)$.
If $x, y\in H$, then $x\ot y$ stands for the operator defined as
$(x\ot y)(z)=\la z,y\ra x$ $(z\in H)$.
By a conjugate-linear isometry we mean a conjugate-linear operator
$U:H\to H$ for which $\| Ux\|=\| x\|$ $(x\in H)$. Clearly, we have
$\la Ux,Uy\ra =\la y,x\ra$ $(x,y\in H)$. Projection means any
self-adjoint idempotent in $B(H)$.

\section{The Results}

Let $\cA\subset B(H)$ be a standard operator algebra, that is, a
*-subalgebra of $B(H)$ which contains every finite rank operator. It is
a folk
result that the automorphisms of $\cA$ are exactly the transformations
of the form
\[
A \longmapsto UAU^*,
\]
where $U\in B(H)$ is a unitary operator.

The main result of the paper reads as follows.

\begin{theorem}
Let $\I \subset B(H)$ be an ideal and $\phi :\I \to \I$ be a
multiplicative
map (linearity or continuity is not assumed). Suppose that for every
$A\in \I$ there exists a sequence $(\phi_n)$
of automorphisms of $\I$ such that $\phi(A)=\lim_n \phi_n(A)$ in the
operator norm.
\begin{itemize}
\item[(1)]
If $\I=B(H)$, then $\phi$ is an automorphism of $\I$.
\item[(2)]
If $\I \subsetneq B(H)$, then there is a (linear) isometry $U$ such that
$\phi(A)=UAU^*$ $(A\in \I)$.
In particular, $\phi$ is a (linear) *-endomorphism of $\I$.
\end{itemize}
\end{theorem}

After this we immediately have the following corollary.

\begin{corollary}
If $F(H)\subsetneq \I \subset B(H)$ is an ideal and $\phi :\I \to \I$ is
a multiplicative map with the property that for every $A\in \I$ there
exists an automorphism $\phi_A$ such that $\phi(A)=\phi_A(A)$, then
$\phi$ is an automorphism of $\I$.
\end{corollary}

\begin{remark}
Let $U\in B(H)$ be an arbitrary isometry. It is easy to see that there
is a sequence $(U_n)$ of unitaries such that $U_n \to U$ in the strong
operator topology. By Banach-Steinhaus theorem it follows that $U_n
AU_n^* \to UAU^*$ in the operator norm for every compact operator $A$.
This shows that the operator $U$ appearing in the statement (2) of our
theorem really only an isometry and not unitary in general.
\end{remark}

\begin{remark}
Clearly, to any isometry $U$ and any finite dimensional subspace
$M\subset H$ there exists a unitary operator $V\in B(H)$ such that
$U_{|M}=V_{|M}$. This implies
easily that for any isometry $U$, the map $A \mapsto UAU^*$ is a local
automorphism of $F(H)$, that is, to every $A\in F(H)$ there is an
automorphism $\phi_A$ of $F(H)$ such that $UAU^*=\phi_A(A)$.
Consequently, if $H$ is infinite dimensional, the
conclusion of the previous corollary does not hold true for $\I=F(H)$.
\end{remark}

\begin{remark}
Observe that since $B(H)$ is a prime ring, from \cite{Mar} we could
get the additivity of any multiplicative map $\phi$ on $B(H)$
if we had supposed that $\phi$ is
surjective. However, as it is seen, this was not the case above.
Clearly, without that surjectivity assumption we do not have
the additivity in general. Consider, for
example, the transformation
\[
A \longmapsto
\left[
\begin{matrix}
I & 0\\
0 & A
\end{matrix}
\right] .
\]
This shows that our results are far from being trivial.
\end{remark}

\begin{proof}[Proof of the Theorem]
Let $\I$ and $\phi$ be as in the statement of the theorem. By the form
of the automorphisms of $\I$ it is obvious that $\phi(P)$ is a rank-one
projection for every rank-one projection $P\in \I$.
We show that $\phi$ is homogeneous on the set of all rank-one operators.
First, let $P$ be a rank-one projection and $\l \in \C$.
Since $\phi(P)$ is also rank-one, we compute
\begin{equation}\label{E:egy}
\phi(\l P)=\phi(P)\phi(\l P)\phi(P)=\mu \phi(P)
\end{equation}
with some $\mu \in \C$.
On the other hand, by the local property of $\phi$, it follows that
there are unitary operators $U_n\in B(H)$ such that
\[
\phi(\l P)=\lim_n U_n (\l P) U_n^* =\l \lim_n U_n P U_n^*,
\]
that is, $\phi(\l P)$ is equal to $\l$ times a rank-one projection.
Comparing this to \eqref{E:egy} we obtain $\l=\mu$ and hence
$\phi(\l P)=\l \phi(P)$.
Now, if $A$ is an arbitrary rank-one operator, then choosing a rank-one
projection $P$ for which $PA=A$ we have
\[
\phi(\l A)=\phi( \l P A)=\phi(\l P)\phi(A)=
\l \phi(P)\phi(A)=
\l \phi(PA)=\l \phi(A).
\]

Let $x \in H$ be any vector. Since, by the local property of $\phi$,
$\phi(x \ot x)$ is a self-adjoint rank-one operator, it follows that
there exists a vector $Tx \in H$ such that  $\phi(x \ot x)=Tx \ot Tx$.
For arbitrary $0\neq x, y \in H$ we have
\[
\phi(x \ot y) = 1/ (\| x\|^2\| y\|^2) \phi(x\ot x \cdot x \ot y \cdot
y \ot y)=
\]
\[
1/ (\| x\|^2\| y\|^2) \phi(x\ot x) \phi(x \ot y) \phi(y \ot y)=
c \cdot Tx \ot Ty,
\]
where $c$ is a complex number (depending on $x,y$).
Since $\phi$ preserves the operator norm (this follows from the local
property of $\phi$), we obtain that
\[
\| x\| \| y\|=\| \phi(x \ot y)\|= |c| \| Tx \ot Ty\| =| c| \| Tx\| \|
Ty\|= |c| \| x\| \| y\|,
\]
that is, $|c|=1$.
Suppose that $\la x, y\ra \neq 0$. We infer
\[
c^2 \cdot Tx\ot Ty \cdot  Tx \ot Ty= \phi(x\ot y )\phi(x \ot y)=
\phi(x \ot y \cdot x \ot y)=
\]
\[
\la x, y \ra \phi(x \ot y)= \la x, y \ra c \cdot Tx \ot Ty
\]
from which it follows that
$c \la Tx, Ty\ra =\la x,y \ra$.
Therefore, $|\la Tx, Ty\ra |=| \la x,y \ra|$.
Similarly, in the case when $\la x , y \ra =0$ we compute
\[
0=\phi(0)=\phi(x \ot y \cdot x \ot y)= c^2 \cdot Tx \ot Ty \cdot Tx \ot
Ty
\]
which gives us that $\la Tx, Ty \ra =0$. Consequently, we have
\[
|\la Tx, Ty\ra |=| \la x,y \ra| \qquad (x,y \in H).
\]
We now apply Wigner's unitary-antiunitary theorem. This celebrated
theorem says that any function on $H$ which preserves the absolute value
of the inner product is of the form $\varphi(x) Ux$ $(x\in H)$, where
$\varphi$ is a so-called phase function (i.e., $\varphi$ is a complex
valued function $H$ of modulus 1) and $U$ is an either linear or
conjugate-linear isometry (see, for example, \cite{Bar}, \cite{Rat},
\cite{ShAl} as well as \cite{Mol} for a new, algebraic approach to the
result). We claim that in our case $U$ is linear.
To see this, suppose on the contrary that $U$ is
conjugate-linear.
Taking into account that $\phi(x\ot x) = Ux \ot Ux$ (recall that
$|\varphi(x)|=1$), in this case we find that
\begin{equation}\label{E:ket}
\la y, x\ra \phi(x\ot y)=\phi(x \ot x \cdot y \ot y)= Ux \ot Ux  \cdot
Uy \ot Uy =
\end{equation}
\[
\la Uy,Ux \ra Ux\ot Uy =\la x,y\ra Ux\ot Uy.
\]
Now, let $x,y, v,w\in H$ be such that $\la x,y\ra , \la x,w\ra, \la v,
w\ra \neq 0$. We compute
\[
\phi(x\ot y \cdot v \ot w)=\la v, y\ra \phi(x \ot w)=
\la v,y\ra \frac{\la x,w\ra}{\la w,x\ra} Ux \ot Uw.
\]
On the other hand, we have
\[
\phi(x\ot y \cdot v \ot w)=\phi(x \ot y) \phi(v \ot w)=
\]
\[
\frac{\la x,y\ra}{\la y,x\ra} Ux\ot Uy
\frac{\la v,w\ra}{\la w,v\ra} Uv\ot Uw =
\frac{\la x,y\ra}{\la y,x\ra}
\frac{\la v,w\ra}{\la w,v\ra} \la y, v\ra Ux\ot Uw .
\]
This implies that
\[
\la x,y\ra \la y, v\ra \la w,x\ra \la v,w\ra=
\la y,x\ra \la v,y\ra \la x,w\ra \la w,v\ra .
\]
Since this obviously does not hold true for every possible choice of the
vectors $x,y, v,w\in H$, it follows that the operator $U$ cannot be
conjugate-linear.
Therefore, we have a linear isometry $U$ on $H$ such that $\phi(x \ot x
)=Ux \ot Ux$ $(x \in  H)$. Similarly as in \eqref{E:ket} we get that
$\phi(x \ot y) =Ux \ot Uy$ if $\la x,y\ra \neq 0$. In the opposite case
choose a unit vector $v \in H$ such that $\la x, v\ra, \la y,v\ra \neq
0$. We have
\[
\phi(x \ot y)=\phi( x\ot v \cdot v \ot y)=Ux \ot Uv \cdot Uv \ot Uy=
Ux \ot Uy.
\]

Now, let $P$ be a projection of rank $n$. Choose pairwise
orthogonal rank-one projections $P_1, \ldots ,P_n$ such that
$P=P_1 + \cdots + P_n$. By the multiplicativity and the local property
of $\phi$ it follows that $\phi(P_1) , \ldots , \phi(P_n)$ are pairwise
orthogonal rank-one projections. Since
\[
\phi(P_i)=\phi(P)\phi(P_i) \phi(P) \leq \phi(P) I\phi(P) =\phi(P)
\qquad (i=1, \ldots, n),
\]
we infer that
\begin{equation}\label{E:har}
\phi(P_1) + \cdots + \phi(P_n) \leq \phi(P).
\end{equation}
By the
equality of the ranks of the operators appearing on the two sides of
\eqref{E:har}, it
follows that we have in fact equality in the above inequality,
that is,
\[
\phi(P_1) + \cdots + \phi(P_n) = \phi(P).
\]
Now, let $A$ be an arbitrary finite rank operator. Let $P$ be a finite
rank projection such that $A=PA$. If $P_1, \ldots ,P_n$ are as above and
$P_1=x_1\ot x_1, \ldots , P_n=x_n \ot x_n$, then we see that
\begin{equation}\label{E:ot}
\phi(A)=\phi(PA)=\phi(P)\phi(A)=
\sum_{i=1}^n \phi(P_i)\phi(A)=
\sum_{i=1}^n \phi(P_i A)=
\end{equation}
\[
\sum_{i=1}^n \phi(x_i \ot A^* x_i)=
\sum_{i=1}^n Ux_i \ot UA^* x_i=
\sum_{i=1}^n U(x_i \ot  x_i) AU^*=
\]
\[
U \bigl(\sum_{i=1}^n x_i \ot  x_i \bigr) A U^*=U(PA)U^*=UAU^*.
\]
If $H$ is finite dimensional, then $U$ is necessarily unitary and hence
we are done in that case. So, in what follows let us suppose that $H$ is
infinite dimensional. From \eqref{E:ot} we see that on $F(H)$, the map
$\phi$ can be represented as
\[ \phi(A)=
\left[
\begin{matrix}
A & 0\\
0 & 0
\end{matrix} \right]
\qquad (A\in F(H)).
\]
Clearly, on the whole ideal $\I$, $\phi$ can
be written as
\[
\phi(A)=
\left[
\begin{matrix}
\psi_{11}(A) & \psi_{12}(A)\\
\psi_{21}(A) & \psi_{22}(A)
\end{matrix} \right]
\qquad (A \in \I).
\]
Let $A\in \I$ and
consider an arbitrary finite rank projection $P$.
From the equality $\phi(AP)=\phi(A)\phi(P)$ we obtain
$\psi_{11}(A)P=AP$, $\psi_{21}(A)P=0$. Since $P$ was arbitrary, it
follows that $\psi_{11}(A)=A$ and $\psi_{21}(A)=0$. Similarly, from the
equality $\phi(PA)=\phi(P)\phi(A)$ we infer that $\psi_{12}(A)=0$.
Consequently, our map $\phi$ is of the form
\begin{equation}\label{E:negy}
\phi(A)=
\left[
\begin{matrix}
 A & 0 \\
 0 & \psi_{22}(A)
\end{matrix} \right]
\qquad (A \in \I),
\end{equation}
where $\psi_{22}$ is obviously multiplicative and it vanishes on the
finite rank operators.

Up till this point $\I$ has been an arbitrary ideal in $B(H)$. Suppose
now that
$\I$ is proper, that is, $\I \subsetneq B(H)$. By the separability of
$H$, the elements of $\I$ are all compact operators.
Let $A$ be an arbitrary compact operator. Denote $s_n(A)$ the $n$th
$s$-number of $A$ which is the $n$th term in the decreasing sequence of
the eigenvalues of the positive compact operator $|A|$, where each
eigenvalue is counted according to its multiplicity.
For a fixed $n \in \N$, denote $\| A\|_n =s_1(A)+\cdots +s_n(A)$.
It is well-known that $\| .\|_n$ is a norm (somtimes called Ky Fan
norm) on $C(H)$ (see \cite[p. 48]{GK}). Therefore, we have
\[
| \|A \|_n-\| B\|_n| \leq \| A-B\|_n \leq n\| A-B\|
\]
for any compact operators $A,B$. By the local property of $\phi$ it
now follows
that $\| \phi(A)\|_n =\| A\|_n$ for every $A\in \I$. This gives us that
the $s$-numbers of $\psi_{22}(A)$ are all zero and hence
$\psi_{22}(A)=0$. Consequently, we have
\[
\phi(A)=
\left[
\begin{matrix}
 A & 0 \\
 0 & 0
\end{matrix} \right]
\qquad (A \in \I)
\]
which yields $\phi(A)=UAU^*$ for every $A\in \I$.

It remains to consider the case when $\I=B(H)$. The function
$\psi_{22}:B(H) \to B(H)$ appearing in \eqref{E:negy} is a
multiplicative map which vanishes on the set
of all finite rank operators. Suppose that $\psi_{22} \neq 0$. It is
easily
seen that $\psi_{22}(P) \neq 0$ for any infinite rank projection $P$.
Indeed, this follows from the fact that for any infinite rank
projection $P$ there is a coisometry $W$ such that $WPW^*=I$.
Choosing uncountably
many infinite rank projections with the property that the product of
any two of them has finite rank
(see, e.g.,  \cite[Proof of Theorem 1]{MolStud}) and taking their values
under $\psi_{22}$ we would
get uncountably many pairwise orthogonal nonzero projections in $B(H)$.
Since this is a contradiction, we obtain $\psi_{22}=0$. So, just as
in the case when $\I \subsetneq B(H)$, we have $\phi(A)=UAU^*$ $(A\in
B(H))$. But due to the local
property of $\phi$ we have $\phi(I)=I$. Therefore, the
isometry $U$ is in fact unitary. This completes the proof of the theorem.
\end{proof}

\begin{proof}[Proof of the Corollary]
It follows from our theorem that there is an isometry
$U$ such that $\phi(A)=UAU^*$ $(A\in \I)$. Since $F(H) \subsetneq \I$,
there exists an operator $A\in \I$
with dense range. By the local property of
$\phi$, $\phi(A)$ must also have dense range. This gives us that $U$ is
surjective.
\end{proof}

To conclude, we feel that it would be an interesting
question to
study our 'multiplicative' reflexivity problem for algebras of
continuous
functions which represent another important class of $C^*$-algebras.

\begin{remark}
The referee of the paper kindly informed us about the article \cite{JL}
where the multiplicative selfmaps of the matrix algebra $M_n(\C)$
are completely determined. He advised us to try to use that result
to reach our conclusion as well as to look for possible generalizations.
Here, we deal with only this second suggestion.

Let us suppose that our multiplicative map
$\phi: \I \to \I$ is continuous in the operator norm topology (this was
not supposed in Theorem).
Assume that for every $A\in \I$ there exists a sequence $(\phi_n)$ of
unstarred automorphisms of $\I$ such that $\phi(A)=\lim_n \phi_n(A)$ in
the operator norm topology. Then similarly to
our theorem we have the following assertions.
\begin{itemize}
\item[(1)]
If $\I=B(H)$, then $\phi$ is an algebra-automorphism of $\I$.
\item[(2)]
If $\I \subsetneq B(H)$, then there are bounded linear operators $T,S$
on $H$ with $ST=I$ such that $\phi(A)=TAS$ $(A\in \I)$.
In particular, $\phi$ is an  algebra-endomorphism of $\I$.
\end{itemize}

We briefly sketch the proof. First we recall that every
algebra-automorphism of $\I$ is of the form $A \mapsto TAT^{-1}$, where
$T$ is an invertible bounded linear operator on $H$. Similarly to the
first part of the proof of our theorem, we can see that $\phi$ is
homogeneous on the set of rank-one operators. One can verify that $\phi$
preserves the rank of the finite rank idempotents and thus we obtain
that for every $n\in \N$,
$\phi$ can be considered as a selfmap of the matrix algebra $M_n(\C)$.
Applying \cite[Theorem 1]{JL} we can infer that $\phi$ is linear on
$F(H)$.
Consequently, $\phi$ is an algebra-endomorphism of $F(H)$. Since $\phi$
maps any operator of rank one into an operator of rank at most one,
and preserves the rank of the idempotents,
following the argument in \cite{Hou} till the proof of Theorem 1.2,
we obtain that $\phi$ is of the form $\phi(A)=TAS$, with some $T,S \in
B(H)$. We should emphasize that this is the place where we need the
continuity of $\phi$. In
fact, it is not too hard to give a rank-preserving algebra-endomorphism
of $F(H)$ which cannot be written in the nice form above. Consider, for
example, a nonsurjective isometry $V$ on $H$ and an arbitrary unbounded
linear operator $L$ whose range is orthogonal to that of $V$. Then
one can check that the map
\[
\sum_{i=1}^n x_i \ot y_i \longmapsto
\sum_{i=1}^n (Vx_i) \ot ((V+L)y_i)
\]
has the desired properties. Turning back to our original problem,
observe that since $\phi$
preserves the idempotents, we have $\la Tx, S^*y\ra =1$ if $\la x,y\ra
=1$. Consequently, $\la Tx, S^*y\ra =\la x,y\ra$ $(x,y \in H)$, that is,
$ST=I$. If
$\I$ is a proper ideal, then $F(H)$ is dense in $\I$. Thus, by the
continuity of $\phi$ we have the
form $\phi(A)=TAS$ on the whole $\I$. It remains to consider the case
when $\I=B(H)$. But this can be treated very similarly to the
corresponding part of the proof of Theorem.

We note that we feel that the complete working-out of the proof above
would not be shorter than what we have seen in the case of our theorem.
\end{remark}

\end{document}